\documentclass[a4paper]{amsart}

\PassOptionsToPackage{usenames,dvipsnames}{xcolor}
\usepackage{multirow,mathabx}

\usepackage{fancyhdr}

\usepackage[utf8]{inputenc}

\usepackage{mathtools}
\usepackage{amsthm,amssymb,latexsym,url,cancel}
\usepackage{pdfpages}

\usepackage{tcolorbox}
\tcbuselibrary{theorems}

\usepackage{tikz}
\usepackage{tikz-qtree}
\usepackage{tkz-euclide}
\usetkzobj{all}
\usetikzlibrary{intersections}
\usetikzlibrary{calc}
\usetikzlibrary{decorations.markings}
\usetikzlibrary{backgrounds,shapes}

\usepackage[all]{xy}
\usepackage{subfig}
\usepackage{enumerate}
\makeatletter
\let\@@enum@org\@@enum@
\def\@@enum@[#1]{\@@enum@org[\normalfont #1]}
\makeatother

\def\co{\colon\thinspace}

\DeclarePairedDelimiter{\form}{\langle}{\rangle}

\newcommand\ba{\begin{align*}}
\newcommand\ea{\end{align*}}
\newcommand\be{\begin{enumerate}}
\newcommand\ee{\end{enumerate}}
\newcommand\bp{\begin{proof}}
\newcommand\ep{\end{proof}}
\newcommand\bpp{\begin{prop}}
\newcommand\epp{\end{prop}}
\newcommand\bpb{\begin{prob}}
\newcommand\epb{\end{prob}}
\newcommand\bd{\begin{defn}}
\newcommand\ed{\end{defn}}
\newcommand\bh{\begin{hint}}
\newcommand\eh{\end{hint}}

\DeclarePairedDelimiter{\pbr}{\{}{\}}

\newcommand\bN{\mathbb{N}}
\newcommand\bR{\mathbb{R}}

\newcommand\bZ{\mathbb{Z}}

\newcommand\CC{\mathcal{C}}

\newcommand\EE{\mathcal{E}}

\newcommand\GG{\mathcal{G}}
\newcommand\PP{\mathcal{P}}

\newcommand\VV{\mathcal{V}}

\newcommand\hull{\operatorname{hull}}

\newcommand\Cay{\operatorname{Cay}}
\newcommand\val{\operatorname{val}}

\newcommand\Isom{\operatorname{Isom}}

\newcommand\Comm{\operatorname{Comm}}

\newcommand{\T}{\mathcal{T}}

\newcommand\sse{\subseteq}

\DeclareMathOperator\Stab{Stab}

\renewcommand{\MR}[1]
{\href{http://www.ams.org/mathscinet-getitem?mr=#1}{MR#1}}

\def\thetitle{{Some groups with planar boundaries}}
\def\theauthors{{Sang-hyun Kim and Genevieve S. Walsh}}
\usepackage{hyperref}
\hypersetup{
   hidelinks = true,
   colorlinks=false,
   plainpages,
   urlcolor=black,
   linkcolor=black
   pdftitle=   \thetitle,
   pdfauthor=  {\theauthors}
}

\theoremstyle{theorem}
\newtheorem{thm}{Theorem}[section]
\newtheorem{lem}[thm]{Lemma}
\newtheorem{cor}[thm]{Corollary}
\newtheorem{prop}[thm]{Proposition}
\newtheorem{con}[thm]{Conjecture}

\newtheorem*{claim*}{Claim}

\theoremstyle{remark}

\newtheorem*{rem}{Remark}
\newtheorem*{hint}{Hint}

\theoremstyle{definition}
\newtheorem{defn}[thm]{Definition}
\newtheorem{prob}{Problem}[section]

\begin{document}
\title[Discussion]\thetitle
\date{\today}

\author[S. Kim]{Sang-hyun Kim}
\address{School of Mathematics, Korea Institute for Advanced Study, Seoul, Korea}
\email{skim.math@gmail.com}
\urladdr{http://cayley.kr}

\author[G. Walsh]{Genevieve Walsh}
\address{Tufts University, USA}
\email{genevieve.walsh@tufts.edu}
\urladdr{}

\begin{abstract}
In this expository note, we illustrate phenomena and conjectures about boundaries of hyperbolic groups by considering the special cases of certain amalgams of hyperbolic groups.
While doing so, we describe fundamental results on hyperbolic groups and their boundaries by Bowditch \cite{Bowditch1998} and Haissinsky \cite{Haissinsky2015IM}, along with special treatments for the boundaries of free groups by  Otal \cite{Otal1992} and Cashen \cite{Cashen2016AGT}.

\end{abstract}

\maketitle

\section{Introduction and Background}

A boundary of a group contains in general, a wealth of information about the group.  For example, if the boundary is a $Z$-set compactification, the dimension of the boundary is one less than the cohomological dimension of the group, and if the group is torsion-free, this boundary is $S^2$ exactly when the group is a $\text{PD}(3)$ group \cite{BM1991}. Perhaps more fundamentally, the boundary can tell us when the group splits over a finite group \cite{Gromovhyp} and when the group splits over a virtually cyclic group \cite{Bowditch1998, Haulmarksplittings}.
In this survey, we will focus on Gromov hyperbolic and relatively hyperbolic groups, and investigate the planarity of their boundaries under specific circumstances. 

 A {\it Gromov hyperbolic group} $G$ is a group that acts {\it geometrically} (co-compactly, properly discontinuously, and by isometries) on some proper hyperbolic space $X$.  The Gromov boundary of such an $X$ is the set of all geodesic rays from a point, where two such rays are equivalent if they have finite Hausdorff distance. 
 The Gromov boundary, $\partial G$, is the boundary of any proper hyperbolic space $X$ that admits a geometric action by $G$.  Crucially, for a hyperbolic group $G$ this boundary is well-defined up to homeomorphism.  Indeed, all such spaces $X$ are quasi-isometric, and since they are hyperbolic, it follows that their boundaries are homeomorphic \cite{Gromovhyp}.  There are several equivalent ways to define the Gromov boundary of a proper hyperbolic metric space; for more extensive details on this definition see Section \ref{sec:prelim} and the excellent survey article \cite{KapBenakli}. 
 
 A {\it relatively hyperbolic group pair} $(G, \mathcal{P})$ is a group which acts \emph{geometrically finitely} on some hyperbolic space $X$, where the set of peripheral subgroups is $\mathcal{P}$;
see Definition \ref{def:relhyp} for a precise definition.
One of the equivalent definitions of a geometrically finite action in the $\Isom(\mathbb{H}^n)$ case  is that $G$ admits a finite-sided fundamental polyhedron \cite{BowditchGF}.  
 A geometrically finite Kleinian group $G$ along with the collection of its maximal parabolic subgroups $\mathcal{P}$ forms a relatively hyperbolic group pair $(G, \mathcal{P})$.  

The \emph{Bowditch boundary} of a relatively hyperbolic group pair $(G, \mathcal{P})$ is the boundary of any proper hyperbolic metric space $X$ such that $(G, \mathcal{P})$ acts on $X$ geometrically finitely. Like a hyperbolic group, for a relatively hyperbolic group {\it pair}
  $(G, \mathcal{P})$ this boundary  
 is well-defined up to homeomorphism. Indeed, the boundary of a relatively hyperbolic group pair coincides with the boundary of the coned-off Cayley graph, with the parabolic fixed points suitably used to compactify this boundary~\cite[Section 7]{Bowditchrelhyp}. 
However, all such spaces $X$ are not quasi-isometric, \cite{Burns}. We will denote this boundary by $\partial_B(G, \mathcal{P})$ to avoid confusion although it is often denoted simply as $\partial(G, \mathcal{P})$. 

 
Much of our intuition about hyperbolic and relatively hyperbolic groups comes from studying geometrically finite Kleinian groups. 
If a topological space embedds in $S^2$ we say that it is \emph{planar}.
The limit set of a geometrically finite Kleinian group is planar.
One may ask if hyperbolic or relatively hyperbolic groups with planar boundaries are virtually Kleinian.  
The Gromov boundaries of one-ended hyperbolic groups do not have cut points.  If one wants to conjecture that every relatively hyperbolic group pair with planar boundary is (even virtually) the fundamental group of a 3-manifold, it is necessary to exclude cut points.   Indeed there are many examples of relatively hyperbolic group pairs with planar Bowditch boundary which are not virtually Kleinian \cite{HWpreprint}.  This can happen when the peripheral subgroups are $\mathbb{Z} \oplus \mathbb{Z}$  by gluing surfaces along their boundaries to a torus, as shown in \cite{HWpreprint}.
 This suggests the following \emph{Planarity Conjectures}: 

\begin{con} \label{conj:hyp} (Cannon, Kapovich-Kleiner, Haissinsky) If a hyperbolic group $G$ has planar Gromov boundary $\partial G$  then $G$ is virtually isomorphic to a Kleinian group. \end{con} 
\begin{con} \cite{HWpreprint} \label{conj:relhyp} If a non-elementary relatively hyperbolic group pair $(G, \mathcal{P})$ has planar Bowditch boundary that does not have cut points, then $G$ is virtually isomorphic to a Kleinian group. \end{con} 

If  $(G, \mathcal{P})$ with $\mathcal{P} = \emptyset$ is a relatively hyperbolic group pair, then $G$ is a hyperbolic group.  So Conjecture \ref{conj:relhyp} is a generalization of the more well-known Conjecture \ref{conj:hyp}. Also see \cite[Corollary 1.4]{GMS} where the authors show that a special and important case of Conjecture \ref{conj:relhyp} is implied by Conjecture \ref{conj:hyp}, namely \cite[Problem 60]{Kapprob}.  Even when a group acts effectively on its boundary and is torsion-free,  it can be virtually Kleinian without being Kleinian ``on the nose".  The first example we know of this phenomena appeared in \cite{KapKleinlow}; more recent examples are in \cite{HSTAJM}.  All of these examples split over cyclic groups, and this is a necessary condition for being both virtually Kleinian and non-Kleinian when the group acts effectively on its boundary; see Proposition \ref{prop:virtual}.
Here we are addressing only a small case of Conjecture \ref{conj:hyp}, for the case of certain amalgams of hyperbolic groups, and limit groups, and this follows readily from known results.  Relatively hyperbolic group boundaries will be used to understand the pieces during the course of the argument. 

\subsection{Plan of the paper} 
In section \ref{sec:examples} we give specific examples of hyperbolic and relatively hyperbolic groups, focusing our examples on the case of geometrically finite Kleinian groups discussed above.  We show how the canonical splittings of these groups can be seen from the Gromov boundary.  In Section \ref{sec:prelim} we give more precise definitions of hyperbolic and relatively hyperbolic groups, and in Section \ref{sec:Bowditch} we give a self-contained synopsis of Bowditch's theory of splittings over $\mathbb{Z}$ for hyperbolic groups.  We also discuss Otal's theorem in Section \ref{sec:Otal} and the relation with relatively hyperbolic group boundaries. It will follow pretty quickly from the results cataloged here that hyperbolic doubles of free groups are virtually Kleinian exactly when their boundaries are planar.

\section{Boundaries of hyperbolic and relatively hyperbolic groups} \label{sec:prelim}
\subsection{Hyperbolic groups}

Let $X$ be a geodesic space, and let $\delta$ be a positive number.
A geodesic triangle $\Delta$ is called \emph{$\delta$--slim} if each side of $\Delta$ is contained in the union of the $\delta$--neighborhoods of the other two remaining sides.
We say such an $X$ is $\delta$-\emph{hyperbolic} if every geodesic triangle in $X$ is $\delta$-slim~\cite{Gromovhyp}.

\begin{defn} A group is \emph{word-hyperbolic} (or simply, \emph{hyperbolic}) if it acts geometrically (properly discontinuously, co-compactly, and by isometries) on a proper $\delta$-hyperbolic metric space for some $\delta>0$.
\end{defn} 

A proper geodesic space is hyperbolic if it is quasi-isometric to a hyperbolic space~\cite[Proposition 2.20]{KapBenakli}.
Using this observation and  Schwarz--Milnor Lemma, the hyperbolicity of a finitely generated group is equivalent to its Cayley graph being $\delta$-hyperbolic for some choice of a finite generating set.

The group boundary $\partial \Gamma$ of a hyperbolic group $\Gamma$  can be defined as
the set of equivalence classes of based geodesic rays where two rays are declared to be equivalent if they have a bounded Hausdorff distance. This boundary is naturally topologized by the basis of the sets of rays that stay close for a long time. See \cite[Section 2]{KapBenakli} for various equivalent definitions of the boundary and its topology.  In particular, we will use the sequential boundary below.  Since quasi-isometric hyperbolic spaces have homeomorphic boundaries  \cite[Proposition 2.20]{KapBenakli} the above definition is equivalent to:

\begin{defn} The boundary $\partial \Gamma$ of a hyperbolic group $\Gamma$ is the topological space $\partial X$ where $X$ is a proper geodesic space on which $\Gamma$ acts geometrically. \end{defn} 

We give a more detailed constuction from the Cayley graph, as $\Gamma$ naturally acts geometrically on its Cayley graph.  \emph{From now on, we let $\Gamma$ be a hyperbolic group, and $\Cay(\Gamma)$ be its Cayley graph with a fixed finite generating set.}

One may extract $\partial\Gamma$ from sequences of group elements.
The \emph{Gromov product} for three points $x,y,z$ in a metric space $(X,d)$ is defined as 
\[
(y,z)_x:=\frac12 \left( d(x,y)+d(x,z)-d(y,z)\right).\]
Then $\partial\Gamma$ coincides topologically with the set of equivalence classes of sequences $\{x_i\}\sse\Cay(\Gamma)$ satisfying \[\lim_{i,j\to\infty} (x_i,x_j)_x=\infty.\]
The meaning of this equation is that the geodesics $[x,x_i]$ stay close to each other for a longer and longer time as $i\to\infty$. Two such sequences $\{x_i\},\{y_j\}$ are equivalent if 
 \[\lim_{i,j\to\infty} (x_i,y_j)_x=\infty.\]
Here, the choice of the base point $x$ is arbitrary and does not alter the topology.
For every pair of distinct points $p,q\in\partial\Gamma$ there exists a geodesic sequence $\{x_i\}\in\Cay(\Gamma)$ such that $\lim_{i\to-\infty} x_i=p$ and $\lim_{i\to\infty} x_i=q$.

By Stallings' theorem, a group $G$ is 
\begin{itemize}
\item \emph{0--ended} if $G$ is finite;
\item \emph{two-ended} if $G$ is virtually infinite cyclic;
\item \emph{$\infty$--ended} if $G$ nontrivially splits as a free product or an HNN-extension over a finite group.
\end{itemize}
The group $G$ is one-ended otherwise.

Recall a topological space is \emph{locally connected} if there exists a basis of open connected sets. 
A \emph{continuum} is a nonempty compact connected metrizable space. In particular, a \emph{Peano continuum} is a continuum that is locally connected.
It is a consequence of deep results by Bestvina--Mess \cite{BM1991} and Bowditch \cite[Corollary 0.3]{Bowconnected} that the boundary of a one-ended hyperbolic group $\Gamma$ is a Peano continuum without a global cutpoint; this latter clause means that $\partial\Gamma\setminus\{x\}$ is connected for all $x\in\partial\Gamma$. 

The \emph{limit set} $\Lambda (H)$ of $H\le\Gamma$ is the smallest nonempty closed $H$--invariant subset of $H$ in $\partial\Gamma$. The set $\Lambda (H)$ can be realized as the set of sequences $p=\{\gamma_n\}\in\partial\Gamma$ such that $\gamma_n\in H$.
A subgroup $H\le\Gamma$ is called \emph{elementary} if it is virtually cyclic, i.e. $0$-- or $2$--ended.

An isometry $f$ of a hyperbolic space $X$ is \emph{loxodromic} if $f$ acts by north-south dynamics on $\partial X$; in particular, $f$ fixes exactly two points on $\partial X$.
It turns out that every infinite order element $g\in\Gamma$ acts {loxodromically} on $\partial \Gamma$. 
If $g\in\Gamma$ is loxodromic, then the limit set $\Lambda(g):=\Lambda(\form{g})$ is a pair of points and $\Gamma$ acts on $\partial\Gamma\setminus\Lambda(g)$ properly cocompactly. 
The limit sets of two loxodromics of $\Gamma$ are either equal or disjoint.

When $\Gamma$ is hyperbolic, we say a subgroup $H\le\Gamma$ is \emph{quasi-convex} if the inclusion $H\hookrightarrow\Gamma$ is a quasi-isometric embedding, see \cite[Corollary 3.6]{BH1999}. 
In this case we have that $H$ is hyperbolic \cite[3.7]{BH1999} and that there exists an $H$--equivariant homeomorphism $\partial H\to\Lambda H$.

If $G,H$ are quasi-convex subgroups of $\Gamma$, then so is $G\cap H$~\cite{Short1990} and moreover, 
\[\Lambda(G\cap H)=\Lambda(G)\cap\Lambda(H).\]

\subsection{Boundaries of 3-manifold groups}
A prime motivational example of a hyperbolic group is a {\it convex co-compact Kleinian group} $G$.  A discrete subgroup of $\Isom(\mathbb{H}^3)$ is called a \emph{Kleinian group}. A Kleinian group is \emph{nonelementary} if it is not virtually cyclic.  

For a Kleinian group $G$, we define its \emph{limit set} $\Lambda(G)$ as the minimal nonempty closed $G$--invariant subset in $\partial\mathbb{H}^3$.
We will denote by $\hull(G)$ the convex hull of $\Lambda(G)$.
We define the \emph{domain of discontinuity} as the complement
\[\Omega_G:=\mathbb{H}^3\setminus\Lambda(G).\] 

A Kleinian group $G$ is called  \emph{convex cocompact Kleinian} if the action of $G$ on $\hull(G)$ in $\mathbb{H}^3$ is cocompact. Equivalently,  the action of $G$ on $\mathbb{H}^3\cup\Omega_G$ is cocompact. 
In the case that $G$ acts freely, 
\[
M_C(G):=\left(\mathbb{H}^3\cup\Omega_G\right)/G\]
is a compact orientable 3-manifold whose boundary consists of hyperbolic surfaces.

When $G$ is a convex co-compact Kleinian group, $G$ acts geometrically on the hyperbolic space $\hull(G)$, so $G$ is a hyperbolic group.   The limit set $\Lambda(G)$ is homeomorphic to $\partial G$, the Gromov boundary of $G$.  

There are some hyperbolic 3-manifolds whose fundamental groups are not hyperbolic as groups.  For example, the fundamental groups of hyperbolic knot complements.  However, these manifolds do admit {\it geometrically finite} hyperbolic structures.  
If a Kleinian group $\Gamma$ admits a finite-sided convex fundamental domain with the side pairings, then we say $M_C(\Gamma)$ is geometrically finite~\cite{BowditchGF}.  
This is equivalent to saying that $M_C(\Gamma)$ is the union of a compact set and a finite number of ``standard cusp regions".  Another notion, which is equivalent for hyperbolic manifold groups $\Gamma$ \cite[Definition GF2]{BowditchGF}, is that every point of the limit set is either a conical limit point or a bounded parabolic point (defined below).  
This is what we will use for the general definition of a relatively hyperbolic group.

\subsection{Relatively hyperbolic groups}

Just as a hyperbolic group is a group which acts geometrically on a hyperbolic metric space, a {\it relatively hyperbolic group} is a group which acts {\it geometrically finitely}  on a hyperbolic metric space. There are many equivalent definitions of a relatively hyperbolic group, see \cite{HruskaRH}. 

\begin{defn} \cite[Definition 1]{Bowditchrelhyp} \label{def:relhyp} Let $G$ be a group, and let $\mathcal{P}$ be a collection of infinite finitely generated subgroups. We say that the group pair $(G, \mathcal{P})$ is a {\it relatively hyperbolic group pair}, or that $G$ is {\it hyperbolic relative to $\mathcal{P}$} if $G$ acts properly discontinuously and by isometries on a proper hyperbolic space $X$ such that 
 
\begin{itemize} 
\item[(i)] Every element of $\partial X$ is either a conical limit point or a bounded parabolic point; 
\item[(ii)] The elements of $\mathcal{P}$ are exactly the maximal parabolic subgroups of $G$. 
\end{itemize} 
In this case, we say that $(G, \mathcal{P})$ acts on $X$ {\it geometrically finitely}. 
\end{defn}

We will shortly define conical limit point, bounded parabolic point, and maximal parabolic subgroup, but we can immediately define the Bowditch boundary, or the relatively hyperbolic boundary of a group pair. 

\begin{defn} The \emph{Bowditch boundary} $\partial_B(G, \mathcal{P})$ 
of a relatively hyperbolic group  $(G, \mathcal{P})$ is  the boundary of any hyperbolic space $X$ that $(G, \mathcal{P})$ acts on geometrically finitely. 
\end{defn} 
It is not true that every such space is quasi-isometric \cite{Burns} but nonetheless, all such $\partial X$ are homeomorphic~\cite[Section 9]{Bowditchrelhyp}.  This defines the boundary of a relatively hyperbolic group pair up to homeomorphism. 

\begin{defn} Suppose a group $G$ acts by homeomorphisms on the boundary $\partial X$ of a hyperbolic space $X$.
\be
\item
A subgroup $P$ of $G$ is {\it parabolic}  if it is infinite, contains no loxodromic, and fixes a point of $\partial X$.   A parabolic subgroup which fixes a point $x_p$ of $\partial X$  is {\it bounded} if $P$ acts cocompactly on $\partial X \setminus \lbrace x_p \rbrace$. 
\item
A point $y\in \partial X$ is a \emph{conical limit point of $G$}
if there exists a sequence $(g_i)_{i \in \mathbb{N}} \subset G$ and distinct points $a,b \in \partial X$ such that $g_i(y) \rightarrow a$,  and such that $g_i(x) \rightarrow b$ for all $x \in \partial X \setminus \lbrace y \rbrace$. \ee
\end{defn}

Let us recall a general definition. A collection of subgroups $\pbr{ H_1,...H_n }$ of a hyperbolic group $G$ is  \emph{almost malnormal} if $H_i \cap gH_jg^{-1}$ is finite for $g\in G$ unless $i=j$ and $g \in H_i$. 
There are two points about the relationship between hyperbolic and relatively hyperbolic structures on the same group that are particularly relevant here. 

\begin{lem} (Bowditch, \cite[Theorem 7.11]{Bowditchrelhyp} If $G$ is a hyperbolic group and $\mathcal{P}$ is a almost malnormal collection of infinite quasi-convex subgroups consisting of finitely many conjugacy classes  then $(G, \mathcal{P})$ is a relatively hyperbolic group pair.  \end{lem} 

\begin{lem} \label{lem:Tranlemma} (Tran, \cite[Main Theorem]{Transquash}) also \cite{Floydsquash} and \cite{Manningsquash}) If $G$ is a hyperbolic group and $(G, \mathcal{P})$ is a relatively hyperbolic group pair, then $\partial_B(G, \mathcal{P})$ is the quotient of $\partial G$ obtained by collapsing each of $\partial P$ for each $P \in \mathcal{P}$ to a point. 
\end{lem}

\section{Examples of hyperbolic and relatively hyperbolic groups and their maximal splittings} 
\label{sec:examples} 
Let us now give some examples of hyperbolic 3-manifold groups, explain their canonical splitting, and illustrate how relatively hyperbolic group boundaries fit into the picture.  Here we will also give examples of Bowditch's canonical splitting (Theorem \ref{thm:bowditch}).  Precise and complete definitions are given in Section \ref{sec:Bowditch}. 

The universal cover of a manifold or a complex, can be a useful tool to determine the Gromov boundary of the fundamental group of the manifold or complex.  As such, we will need to deal with the components of the pre-image of certain curves. 

\begin{defn} Let $W \subset M$.  If $\widetilde M \rightarrow M$ is a cover, then the components of the pre-image of $W$ in $\widetilde M$ are called the {\it elevations} of $W$. 
\end{defn} 

\subsection{Example 1: Three surfaces glued along a circle} 
Let $S_{g,1}$ denote a hyperbolic surface of genus $g\ge1$ with one boundary component. 
We glue $S_{1,1}, S_{2,1}$ and $S_{3,1}$ along their boundaries to one copy of $S^1$ (call it $a$) by degree 1 maps.
Call the resulting 2-complex $\T$ (for Triple).   This 2-complex $\T$ embeds into $\mathbb{R}^3$, and one can take a regular neighborhood in $\mathbb{R}^3$ to obtain a 3-manifold $M$ with boundary.  Since the fundamental group $\Gamma$ of $\T$ acts geometrically on a hyperbolic space (the universal cover $\widetilde \T$ of the 2-complex $\T$) the group $\Gamma$ is hyperbolic.
Therefore $M$ is realized as a hyperbolic 3-manifold with boundary,
and $\Gamma$ is a convex cocompact Kleinian group~\cite[Theorem 2.24]{seriescrash}.
The resulting 3-manifold $M$ cannot be realized as a hyperbolic 3-manifold with {\it totally geodesic} boundary, as there are many essential annuli.  
Both the convex hull of the limit set of $\Gamma$ and $\widetilde \T$ are proper hyperbolic metric spaces on which $\Gamma$ acts geometrically, and so,
 $\Lambda(\Gamma) \cong \partial \widetilde \T \cong \partial \Gamma$. 

The subgroups corresponding to the fundamental groups $G_1$, $G_2$ and $G_3$ of the surfaces $S_{1,1}$, $S_{2,1}$ and $S_{3,1}$ are free quasi-convex subgroups.  Pick an elevation of $a$ and let its stabilizer be generated by an element $g_a \in \Gamma$. 

We consider two representations of the groups $G_i$ as a discrete group of isometries of $\mathbb{H}^2$, as this will be helpful in visualizing the universal cover $\widetilde \T$.  The Fuchsian group $G_1$ acts on $\mathbb{H}^2$, and acts co-compactly on a convex subset $C$ of $\mathbb{H}^2$ where the boundary curves of $C$ correspond to elevations of boundary curve of $S_{1,1}$. We continue to denote the stabilizer of one of the elevations in $\mathbb{H}^2$ by $\form{g_a}$. The conjugates of the cyclic subgroup $\form{g_a}$ act loxodromically, and each leaves an elevation of $a$ invariant. In this representation, the limit set is a Cantor set. It is naturally a subset of $S_1$ and we say that this Cantor set is {\it cyclically ordered}; see Proposition~\ref{p:sep}. Also, $G_{1}$ acts on $\mathbb{H}^2$ as a finite co-volume Fuchsian group, where conjugates of $\form{g_a}$ act parabolically. In this representation, the limit set is a circle.  The way to get from the limit set of the first  representation to the limit set of the second representation is to collapse the endpoints of the loxodromic boundary elements to points, just as in Lemma \ref{lem:Tranlemma}. The same process can be done with the other surface groups. 

Now consider $\widetilde \T$ as a union of convex pieces of hyperbolic planes (each the universal cover of one of the $S_{i,1}$) glued together along their boundary curves.  Any elevation of one of the $S_{i,1}$ is isometric to a convex subset of $\mathbb{H}^2$, and is a convex subset of $\widetilde \T$. Its boundary is the limit set of a conjugate of $G_i$, and is a Cantor set.  This is an example of a cyclically ordered Cantor set.  When we collapse the pairs of points of this Cantor set that correspond to the endpoints of conjugates of $g_a$ in $S_{i,1}$ we get a circle, and this is the largest such group with this property. This is called a {\it maximal hanging Fuchsian group}.   When we are looking at $S_{2,1}$ or $S_{3,1}$, we could have taken a subgroup of one of these groups corresponding to a lower genus subsurface. Then the limit set of the subgroup would be a cyclically ordered Cantor set, and when we collapsed the boundaries of the subgroup-invariant collection of cyclic subgroups (corresponding to the boundary of this subsurface) we would obtain a circle. However, it would not be a maximal such subgroup. See the formal definition of a maximally hanging Fuchsian subgroup, Definition \ref{def:maxhang} in Section \ref{sec:Bowditch}. 

Apart from maximal hanging Fuchsian subgroups, another important collection of subgroups is the collection of cyclic subgroups which are conjugates of $\form{g_a}$. Each of these groups stabilizes a line in the 2-complex $\widetilde \T$, an elevation of $a$.   The lines are where the pieces of hyperbolic planes are glued together in $\widetilde \T$.  Removing one of these elevations will break the complex $\widetilde \T$ into three pieces.  The two points on $\partial \widetilde \T$ which are the boundary of one of these elevations will break the boundary into three components.  They will be a pair of local cut points (each of valence 3) that forms a cut pair, which will separate the boundary into three pieces; see Lemma~\ref{l:lc-gamma}.

The complex  $\widetilde \T$ consists of elevations of the $S_{i,1}$'s, glued along elevations of $a$.  We can embed a bipartite tree in $\widetilde \T$, by putting a black vertex in each elevation of the $S_{i,1}$, and a white vertex in each elevation of $a$.  Here the white vertices correspond to elementary subgroups of this splitting, and the black vertices correspond to the non-elementary vertices (Definition \ref{d:2ended}).  Then we connect the associated vertices whenever an elevation of $a$ meets an elevation of $\widetilde \T$.  The tree is exactly the tree corresponding to the maximal splitting tree in Section \ref{sec:Bowditch}. The associated graph of groups will be a tripod, where the center vertex is the cyclic group corresponding to $\form{g_a}$, the edge groups are all $\mathbb{Z}$, and the hanging vertex groups are the maximal hanging Fuchsian vertex groups described above. 

Now $\widetilde \T$ is not canonical for the hyperbolic group $\Gamma$, but the boundary is. From this boundary, one can read a tree; see Theorem \ref{thm:bowditch}. The boundary of  $\Gamma$ consists of the boundaries of the elevations of the $S_{i.1}$, glued together along the boundaries of the elevations of $a$, compactified by the endpoints of this tree.   To see this last fact note that every sequence of points tending to infinity in $\Gamma$ goes through a sequence of elevations of $a$ and elevations of the $S_{i,1}$.  This sequence either terminates, in which case the sequence is associated to the boundary of some subgroup labeled by a white or black vertex, or it does not terminate, in which case it is associated to the boundary of the bipartite tree.

\subsection{Example 2: A hyperbolic manifold with totally geodesic boundary glued to a surface} Suppose $M$ is the union of two hyperbolic 3-manifolds glued along an essential annulus as follows.
Let $M_1$ be a hyperbolic manifold with a totally geodesic boundary,
  and let  $M_2 = S_{g,1} \times I$.
 We glue these together by gluing the annulus $\partial S_{g,1} \times I$ to a neighborhood of an essential simple closed curve $b$ on $\partial M_1$, and denote the resulting 3-manifold as $M$.
By the Bestvina--Feign combination theorem~\cite{BF1992},  the resulting group $\Gamma=\pi_1(M)$ is hyperbolic and $M$ is atoroidal.
Since $\partial M \neq \emptyset$, we may apply Thurston's hyperbolization of Haken manifolds.  By hyperbolization, $M$ can be realized as a hyperbolic 3-manifold \cite[Theorem 2.24]{seriescrash}, \cite[Theorem A']{Morgan} and $\Gamma$ can be realized as a convex cocompact Kleinian group. 
The group $\Gamma$ admits a graph of groups decomposition $\pi_1(M_1) *_\mathbb{Z} \pi_1(S_{g,1}) $, and this splitting is visible from the limit set $\Lambda(\Gamma)$, which is homeomorphic to the the Gromov boundary $\partial \Gamma$.

Instead of building a 3-manifold, we can also build a negatively curved complex whose fundamental group is $\Gamma$. 
Let $\pi_1(S_{g,1})$ act geometrically on a convex subset $C_2$ of $\mathbb{H}^2$ as the surface groups do in Example 1.
Recall $\pi_1(M_1)$ acts geometrically on $C_1$, the convex hull of its limit set in $\mathbb{H}^3$.  
We will require that the length of the curve $b$ in $S_{g,1}$ is the same as the length of the curve $b$ on the boundary of $M_1$.  Let $D$ be the complex obtained by gluing $S_{g,1}$, 
realized as a geometric quotient of the corresponding convex hull $C_2$, to $M_1$, realized as a geometric quotient of $C_1$.  Then the universal cover $\tilde D$ is a union of copies of $C_1$ and $C_2$ (elevations of $M_1$ and $M_2$), glued together along elevations  of $b$.  
Call the stabilizer of one the elevations of $b$ as $\gamma_b$, and note that the rest of the stabilizers are exactly the conjugates of $\form{\gamma_b}$.   The boundary of any copy of $C_2$ is a cyclically ordered Cantor set, such that when one collapses the endpoints of the quasi-convex groups corresponding to conjugates of $\form{\gamma_b}$, the result is a circle, and these are maximal for the property.  The end points of the conjugates of $\form{\gamma_b}$ are local cut points of valence two, and each pair cuts the boundary of $\tilde D$ into two pieces as the associated elevation of $b$ cuts the complex $\tilde D$ into two pieces. This splitting has three types of pieces, the cyclic groups corresponding to the conjugates of $\form{\gamma_b}$, the maximally hanging Fuchsian groups corresponding to the boundaries of the copies of $C_2$, and the {\it rigid} pieces corresponding to the boundaries of the copies of $C_1$.  

Again we can build the bipartite splitting tree of Section \ref{sec:Bowditch}  which embeds in the complex $\widetilde D$ by putting a black vertex in each elevation of $M_1$ and of $M_2$.  These will be the non-elementary vertices in Definition \ref{d:2ended}.  We also put a white vertex in each elevation of $b$, and these are the elementary vertices.  Again we connect the vertices when the associated elevations meet in the universal cover.  The vertex stabilizers are the same as the stabilizers of the associated elevations.  In this case there are two types of non-elementary vertices: namely, the rigid vertices associated to the hyperbolic 3-manifold with totally geodesic boundary and the hanging Fuchsian vertices. As in Example 1, every point of the boundary of the hyperbolic group is either in the limit set of one of the stabilizers of a vertex in the canonical tree of Bowditch in Theorem \ref{thm:bowditch} or can be associated with an endpoint of this tree. 

In this example, the stabilizers of each elevation of $M_1$ (which is isometric to the convex hull of a Sierpinski carpet) is a {\it rigid} subgroup - it does not split over any two-ended group. Note that $M_1$ is also rigid in the sense that there is only one hyperbolic structure with totally geodesic boundary.  One can see this by doubling $M_1$ over the totally geodesic boundary, obtaining a closed hyperbolic manifold, and applying Mostow rigidity. More relevant to the work here, is that there is one hyperbolic structure on an elevation $C_1$ of $M_1$ where the elevations of $b$ in $C_1$ are parabolic.

\begin{defn} \label{d:rigid} We say that a group $G$ {\it  splits with respect to a collection of subgroups $\mathcal{A}$}  if $G$ admits a graph of groups decomposition where every subgroup in $\mathcal{A}$ is contained in some vertex group of the associated Bass-Serre tree. Equivalently, every subgroup in $\mathcal{A}$ is conjugate into one of the vertex groups of the graph of groups.  Let $G$ be a hyperbolic group with connected Gromov boundary $\partial G$, and $\mathcal{A}$ a collection of two-ended subgroups.  We say that a subgroup $R$ of a hyperbolic group $G$ is {\it rigid with respect to $\mathcal{A}$}  if it does not split over two-ended subgroups with respect to $\mathcal{A}_R = \lbrace A \cap R: A \in \mathcal{A} \rbrace$.   
\end{defn} 

Compare the above with Definition \ref{d:rigidrel}. The rigid subgroups in our example here do not split over any two-ended subgroup; we will see in the next example that there are rigid subgroups which can split in many different ways, for example a free group.  However, they do not split over a 2-ended group with respect to the collection $\mathcal{A}$ of prescribed virtually cyclic groups. 


\subsection{Example 3: The double of a free group}\label{ex:double}
Let $F$ denote a nonabelian free group written as 
\[F:=\form{a_1,\ldots,a_n}.\]
We fix a copy $\bar F$ of $F$, and let $\sigma\co F\to \bar F$ be an isomorphism. 
For $g\in F$, we also write $\sigma(g)=\bar g$.
For each word $w\in F$, we define the \emph{(Baumslag) double of $F$ along $w$} as the group
\[D(F,w)=\form{
a_1,\bar a_1,\ldots, a_n,\bar a_n, \mid
w=\bar w}.\]

A word $w\in F$ is \emph{root-free} if it is not a proper power of another word. 
The doubles are among the earliest examples of hyperbolic groups; namely, the double $D(F,w)$ is hyperbolic if and only if $w$ is root-free.
One can further characterize group theoretic properties of $D(F,w)$ using the following terms.

\begin{defn}\label{defn:acyl}
Let $w\in F$ be a word. 
\be
\item
We say $w$ is \emph{indecomposable} (or, \emph{diskbusting}~\cite{Canary1993}) if there does not exist a nontrivial free product decomposition $F=A\ast B$ such that  $w$ is conjugate into $A$ or $B$.
\item
We say $ w$ is \emph{acylindrical} if one cannot write $F$ as the free product decomposition $F=A\ast_C B$ or an amalgamated HNN decomposition $F=A\ast_C$ such that  $C$ is 2-ended
and such that $w$ is conjugate into $A$ or $B$.
\ee
\end{defn}
Note that $w$ is non-acylindrical if $F$ admits a nontrivial graph of groups decomposition with two-ended edge groups such that $w$ is conjugate into one of the vertex groups.  In other words, the word $w$ is acylindrical exactly when $F$ is rigid with respect $\form{w}$.
It is well-known that $D(F,w)$ is one-ended if and only if $w$ is indecomposable (cf. \cite{GW2010,KW2012}]).

We will identify $F$ with the fundamental group of a genus $n$ handlebody $H$.
A word $w\in F$ is called \emph{geometric} if $w$ can be realized by a simple closed curve on $\partial H$.  

Let us now assume $w$ is indecomposable, acylindrical, geometric and root-free in $F$, 
realized by a simple closed curve $\gamma\sse \partial H$.
In particular, $D(F,  w)$ is a one-ended hyperbolic group. We will see how acylindricity comes into play in a moment.
We let $H_1$ and $H_2$ be two copies of $H$, and glue $H_1$ and $H_2$ along an annulus $A_\gamma$ on each of $\partial H_i$ with $\gamma$ as its core.
Denote by $M_\gamma$ the 3-manifold thus obtained.  The properties of the manifold will  depend on properties of the word that represents $\gamma$ up to conjugacy in the free group $\pi_1(H)$.
 
The group $D(F,w)$ is hyperbolic. By hyperbolization (as in Example 2 above, $\partial M \neq \emptyset$ and $\pi_1(M)$ is a hyperbolic group)  $M_\gamma$ is a hyperbolic manifold, the union of $H_1$ and $H_2$ glued along the annulus $A_\gamma$.  We can realize $\pi_1(M_\gamma)= D(F,w)$ as a convex cocompact Kleinian group such that the limit set $\Lambda(D(F,w))$ is connected, as $D(F,w)$ is one-ended.


Consider the universal cover $\widetilde M_\gamma$ of $M_\gamma$. This can be realized as the convex hull of $\Lambda(D(F,w))$.  It consists of elevations of $H_1$ and $H_2$, glued along the elevations of $A_\gamma$.  Each elevation of $A_\gamma$ has two points, and these are the fixed points of the infinite cyclic group stabilizing this elevation.  Each such pair of endpoints separates $\partial\widetilde M_\gamma = \Lambda(D(F,w))$ into two pieces, and each elevation of $A_\gamma$ separates $\widetilde M_\gamma$ into two pieces. We can form the canonical splitting tree of Theorem \ref{thm:bowditch} by putting a black vertex in each elevation of $H_1$ and $H_2$, and putting a white vertex in each elevation of $A_\gamma$.  Then we connect vertices when the associated elevations meet.  We claim that all stabilizers of the black vertices are rigid with respect to the stabilizers of white vertices that are incident, as in Definition \ref{d:rigid}.  Indeed, each stabilizers of an elevation of $H_1$ = $G_v$ is a free group $F_n$, and the stabilizers of the incident vertices are the conjugates of $\form{w}$ in this $F_n$. Because the word $w$ is acylindrical, this $F_n$ does not split over a two ended-group where $\form{w}$ is conjugate into one of the vertex groups. So each black vertex stabilizer is rigid with respect to the stabilizers of the white vertex groups which are incident to it. 

The stabilizers of the elevations of the $H_i$ are rigid in the hyperbolic manifold sense as follows.  Let $\mathcal{H}$ be an elevation of one of the $H_i$.  The quotient of $\mathcal{H}$ by its stabilizer is a handlebody.  Let us make the stabilizers of the elevations of the annulus $A_\gamma$ parabolic, and let $M_H = \mathcal{H}/\Stab(\mathcal{H})$ be the resulting hyperbolic manifold with totally geodesic boundary. 
There are no essential annuli in $M_H$; indeed, if there were essential annuli in $M_H$ then there would be a splitting with $ w$ conjugate into some vertex group, violating the acylindricity. Then we have the uniqueness of the hyperbolic structure with totally geodesic boundary roughly as follows. Consider the manifold $M_H = \mathcal{H}/\Stab(\mathcal{H})$, where the conjugates of $\form{w}$ are parabolic. Since $w$ can be realized as a curve $\gamma$ on the boundary of $H_1$, the manifold $M_H$ has boundary component(s) which are $\partial H_1$ with the curve $\gamma$ parabolic. This will be a cusp on the boundary.  The double of the hyperbolic manifold $M_H$ along its cusped boundary does not contain any essential tori and so admits one (up to conjugation in $\text{PSL}(2, \mathbb{C})$)  complete finite volume hyperbolic structure by Mostow-Prasad rigidity.  This manifold admits an isometry which fixes the boundary of $M_H$.  Thus the manifold $\mathcal{H}/\Stab(\mathcal{H})$ admits a  unique hyperbolic structure with totally geodesic boundary.  Note that the Bowditch boundary of the relatively hyperbolic group pair $(F_n, \form{w})$ is the limit set of $\Stab(\mathcal{H})$ where the $w$ conjugates are parabolic.  By Tran's Lemma \ref{lem:Tranlemma} this is obtained as a quotient of the Cantor set by pinching the endpoints of the conjugates of $w$.

\section{Bowditch's canonical splitting of hyperbolic groups}  \label{sec:Bowditch}
One of the most important tools analyzing boundaries of hyperbolic groups is Bowditch's \emph{canonical splitting}. Let us exhibit a self-contained definition of this splitting and summarize its key algebraic features, following~\cite{Bowditch1998}.


\subsection{Elementary splittings}
A splitting of a group is often considered as a finite graph of groups decomposition~\cite{Scott1979}. 
However, an equivalent definition using an action on a tree seems more apt when one compares various splittings of a given group~\cite[Section 6]{Bowditch1998}. 


Recall an action (simplicial, by default) of a group $\Gamma$ on a tree $\Sigma$ is \emph{minimal} if $\Sigma$ has no proper nonempty $\Gamma$--invariant subtree.
The action is \emph{co-finite} if the quotient $\Sigma/\Gamma$  is finite. 
A graph $\Sigma=(\VV,\EE)$ is \emph{bipartite} on $(X,Y)$ if the vertex set has a partition $\VV=X\coprod Y$ and the edge set $\EE$ is a symmetric binary relation satisfying
\[
\EE\sse (X\times Y)\cup (Y\times X).\]
We simply express this by saying $\Sigma=(X\coprod Y,\EE)$ is bipartite.
If $\Gamma$ acts on $\Sigma$, then the stabilizer group of a vertex or an edge $x$ is written as $\Gamma(x)$.

A subgroup $H$ of a nontrivial hyperbolic group is \emph{elementary}  if $H$ is 0-or two-ended.
In particular, $H$ is  \emph{maximal elementary} if it is maximal among elementary subgroups.
\bd\label{d:2ended}
Let $\Gamma$ be a one-ended group.
By an \emph{elementary splitting} of  $\Gamma$, we mean a minimal, co-finite action of $\Gamma$ on a simplicial bipartite tree \[\Sigma=(\VV=\VV_{\mathrm{e}}\coprod\VV_{\mathrm{ne}},\EE)\] such that the following hold.
\be[(i)]
\item\label{p:vertex} (vertices) 
Distinct vertices have distinct stabilizer groups.
\item\label{p:bipartite} (bipartite)
We have that \begin{align*}
\VV_{\mathrm{e}}&=\{v\in\VV\mid \Gamma(v)\text{ is elementary}\},\\ 
\VV_{\mathrm{ne}}&=\{v\in\VV\mid \Gamma(v)\text{ is nonelementary}\}.
\end{align*}
\ee
Furthermore, if $\Gamma(v)$ is maximal elementary for each $v\in\VV_{\mathrm{e}}$, then 
we say $\Sigma$ is a \emph{maximal elementary splitting}.
\ed

We remark that the tree is allowed to be locally infinite.
The term \emph{elementary splitting} refers to the condition that each edge stabilizer group is necessarily elementary.
As $\Gamma$ is one-ended, each $\Gamma(v)$ is two-ended (i.e. not 0-ended) for each $v\in\VV_{\mathrm{e}}$ by Stallings' Theorem,

The commensurator group of a subgroup $H\le\Gamma$ is defined as 
\[
\Comm_\Gamma(H)=\{ g\in\Gamma\mid [H:H\cap H^g]<\infty\text{ and }[H^g:H\cap H^g]<\infty\}.\]
We say $H$ is \emph{full} if $\Comm_\Gamma(H)=H$.
If $H$ is quasi-convex in $\Gamma$, then we have that 
\[
\Comm_\Gamma(H)=\{g\in\Gamma\mid g\Lambda (H) =\Lambda (H)\},\]
which coincides with the unique maximal finite-index extension of $H$ in $\Gamma$.
Moreover, $\Comm_\Gamma(H)$ is full and quasi-convex.
If $g\in\Gamma$ is loxodromic, then 
\[
\Comm_\Gamma(g):=\Comm_\Gamma(\form{g}) = \{h\in\Gamma\mid h g^n h^{-1}=g^{\pm n}\text{ for some }n\in\bN\}.\]
It follows that $\Comm_\Gamma(g)$ is the unique maximal elementary subgroup of $\Gamma$ containing the loxodromic element $g$. 

\begin{rem}
More generally, if $\Gamma$ is hyperbolic but not assumed to be one-ended, then it is reasonable to define an elementary splitting  $\Sigma$ of $\Gamma$ after replacing the condition (i) above by the following.
\be
\item[(i)'] 
If two distinct vertices have the same stabilizer group $H$, then $H$ is finite.
\ee
If we further assume that each vertex stabilizer group of $\Sigma$ is full and quasi-convex, then the above description of commensurator groups implies that the limit sets of distinct vertex stabilizer groups will be distinct unless those limit sets are empty.
\end{rem}
An elementary splitting yields a usual finite graph of groups decomposition from quotienting by $\Gamma$; see~\cite{Serre-trees} for instance. 
Conversely, if a one-ended hyperbolic group is written as a graph of groups such that its induced Bass--Serre tree action satisfies the condition  (\ref{p:bipartite})  above, then
 it is often easy to enforce the condition (\ref{p:vertex}) after consolidating certain vertices as follows.
\begin{prop}\label{p:consol}
Let $\Gamma$ be a one-ended hyperbolic group admitting an action on a simplicial bipartite tree 
\[\Sigma=(\VV_{\mathrm{e}}\coprod\VV_{\mathrm{ne}},\EE)\] 
such that the condition (ii) above holds and such that each vertex stabilizer group is full.
Then there exists another tree $\Sigma'$ and a surjective graph map
\[
\Sigma\rightarrow\Sigma'\]
such that the induced action 
of $\Gamma$ on $\Sigma'$ satisfies (i) and (ii). 
Furthermore, if a vertex $v$ maps to a vertex $v'$ by this graph map then we have
\[\Gamma(v)=\Gamma(v').\]
\end{prop}
\bp
As we stated above, we may consider the condition (i)' instead of (i).
Let us first consider two distinct vertices $u$ and $v$ having the same stabilizer groups. 
If $u\in\VV_{\mathrm{ne}}$, then the open geodesic interval $(u,v)$ contains a vertex $w\in\VV_{\mathrm{e}}$. It would then follow that 
\[\Gamma(u)=\Gamma(u)\cap\Gamma(v)\le\Gamma(w),\]
which is a contradiction. 
Hence we have that $u,v\in\VV_{\mathrm{e}}$. Furthermore,  the same reasoning reveals that the stabilizer group of each vertex in $\VV_{\mathrm{e}}\cap [u,v]$ coincides with $\Gamma(u)=\Gamma(v)$.

We also need the following observation. Suppose two vertices $u\in\VV_{\mathrm{e}}$ and $v\in\VV_{\mathrm{ne}}$ satisfy that
\[
H:=\Gamma(u)\cap\Gamma(v)\]
is infinite.
Let $w\in \VV_{\mathrm{e}}$ be the neighbor of $v$ in the geodesic interval $[u,v]$;
possibly, we have $u=w$. Since the infinite elementary group $H$ is contained  in both of the maximal elementary groups $\Gamma(u)$ and $\Gamma(w)$, it follows from fullness, which implies uniqueness, that  $\Gamma(u)=\Gamma(w)$.

Let us identify the vertices of $\Sigma$ having the same stabilizer groups, and obtain
\[
\VV_{\mathrm{e}}':=\VV_{\mathrm{e}}/\!\!\sim,\quad
\VV_{\mathrm{ne}}' :=\VV_{\mathrm{ne}}/\!\!\sim.\]
Moreover, we declare that $[u]\in \VV_{\mathrm{e}}'$ and $[v]\in \VV_{\mathrm{ne}}'$ are adjacent if and only if $\Gamma(u)\cap\Gamma(v)$ is infinite. We let 
\[\Sigma'=(\VV_{\mathrm{e}}'\coprod\VV_{\mathrm{ne}}',\EE')\]
denote the resulting graph.

We see that two vertices $[x]\in \VV_{\mathrm{e}}'$ and $[y]\in \VV_{\mathrm{ne}}'$ are adjacent in $\Sigma'$
if and only if some representative of $[x]$ is adjacent to $y$ in $\Sigma$.
Moreover, adjacent vertices in $\Sigma$ map to adjacent vertices in $\Sigma'$ by the natural quotient map.
We have an induced action of $\Gamma$ on $\Sigma'$, satisfying $g.[u]=[g.u]$ for each $g\in \Gamma$. 

Let $v$ be a vertex of $\Sigma$. If $g\in \Gamma$ fixes $[v]$ in $\Sigma'$, then by definition we have 
\[ g\Gamma(v)g^{-1}=\Gamma(g.v)=\Gamma(v).\]
The fullness of $\Gamma(v)$ implies that $g\in\Gamma(v)$ and hence,
\[\Gamma([v])=\Gamma(v).\]

It remains to prove that $\Sigma'$ is indeed a tree.
Assume for contradiction that there exists a nontrivial reduced cycle $C'$ in $\Sigma'$. By the previous  paragraphs, $C'$ can be represented as the union of a sequence of length-two paths in $\Sigma$ as follows:
\[ (u_0, u_1,u_2), (u_2',u_3,u_4), (u_4',u_5,u_6),\ldots,(u_{2k-2}',u_{2k-1},u_0').\]
Here, we have that $u_{2i}\in\VV_{\mathrm{e}}$ and that 
\[\Gamma(u_{2i})=\Gamma(u_{2i}').\]
We connect $u_{2i}$ and $u_{2i}'$ by a geodesic $\alpha_i$ in $\Sigma$, and assume that the total length of the resulting cycle $C$ in $\Sigma$ is minimal.

Since $\Sigma$ is a tree, the cycle $C$ has a backtrack. Without loss of generality, we may assume that the backtrack occurs at $u_2$, and that $\alpha_2=(u_2,u_1,v,w,\ldots,u_2')$.
Then we could have reduced the total length of $C$ by replacing
\[ (u_0, u_1,u_2),\quad  \alpha_2=(u_2,u_1,v,w,\ldots,u_2'),\quad (u_2',u_3,u_4)\]
by
\[ (u_0, u_1,v), (v,w,\ldots,u_2'),(u_2',u_3,u_4).\]
This contradicts the minimality, and we conclude that $\Sigma'$ is  a tree. 
\ep

We now list key combinatorial properties of a maximal elementary splitting, motivated by the above proposition.

\begin{prop}\label{prop:basic}
For a maximal elementary splitting $\Sigma=(\VV_{\mathrm{e}}\coprod\VV_{\mathrm{ne}},\EE)$ of a one-ended nonelementary hyperbolic group $\Gamma$, the following hold.
\be
\item\label{p:prop-edge}
Each edge stabilizer group is two-ended and quasi-convex.
\item\label{p:prop-vertex}
Each vertex stabilizer group is full and quasi-convex.
\item\label{p:prop-adjacent}
Let $u,v$ be distinct vertices. 
Then  $\Gamma(u)\cap \Gamma(v)$ is infinite if and only if either $u$ and $v$ are adjacent, or $u$ and $v$ have a common elementary neighbor; in this case,  $\Gamma(u)\cap \Gamma(v)$ is two--ended.
\item\label{p:prop-incident}
Let $e,f\in\EE$. Then $\Gamma(e)\cap\Gamma(f)$ is two-ended if and only if $e$ and $f$ share a vertex from $\VV_{\mathrm{e}}$.
\ee
\end{prop}
\bp
(\ref{p:prop-edge}) This follows from Stallings' theorem on the number of ends.

(\ref{p:prop-vertex}) Since each edge group is quasi-convex, so is each vertex group; see~\cite[Proposition 1.2]{Bowditch1998}. 
In order to see the fullness, we let $v\in\VV_{\mathrm{ne}}$ and $H:=\Gamma(v)$. 
For each $g\in \Comm_\Gamma(H)$, it suffices to show that $g\in H$.

Assume $g\not\in H=\Gamma(v)$, so that  $g.v\ne v$.
Then the geodesic interval $[v,g.v]$ in $\Sigma$ contains some $u\in\VV_{\mathrm{e}}$
such that $H\cap gHg^{-1}\le \Gamma(u)$. This contradicts the assumption that $\Gamma(u)$ is elementary
and that $g$ is a commensurator of $H$.



(\ref{p:prop-adjacent})  
Suppose $\Gamma(u)\cap\Gamma(v)$ is infinite. Since there exists at most one maximal elementary subgroup of $\Gamma$ containing $\Gamma(u)\cap\Gamma(v)$, we see that $u$ and $v$ cannot be both elementary. So, we may assume $v\in\VV_{\mathrm{ne}}$. Then the closed interval $[u,v]\sse \Sigma$ contains some $w\in \VV_{\mathrm{e}}$ that is adjacent to $v$. 
If $u=w$, then we are done. So, we suppose $u\ne w$. Since $\Gamma(u)\cap\Gamma(v)\le\Gamma(w)$, the uniqueness argument again implies that $u\in\VV_{\mathrm{ne}}$
and that the interval $(u,w)\sse\Sigma$ does not contain elementary vertices. It follows that $w$ is a common elementary neighbor of $u$ and $v$. 

To see the converse, it suffices to consider the case that $u$ and $v$ have a common elementary neighbor $w$. Then $\Gamma(u)$ and $\Gamma(v)$ both contain some two-ended subgroups (namely, corresponding edge groups) of the maximal elementary group $\Gamma(w)$. Since these two-ended subgroups are commensurable, it follows that $\Gamma(u)\cap\Gamma(v)$ is infinite. Furthermore, $\Gamma(u)\cap\Gamma(v)$ is two-ended as it is contained in $\Gamma(w)$. 



The proof of (\ref{p:prop-incident}) is very similar.
\ep

\begin{rem}\label{rem:determine}
By part~(\ref{p:prop-adjacent}) above, we see that the graph structure of  $\Sigma$ is uniquely determined by the set of its vertex stabilizer groups. \end{rem}

\subsection{Canonical splitting of $\Gamma$}
\emph{From now on, let us assume that $\Gamma$ is a one-ended hyperbolic group that is not a cocompact Fuchsian group.}
Let us describe Bowditch's \emph{canonical splitting} $\Sigma=(\VV,\EE)$ through the action of $\Gamma$ on $\partial \Gamma$.

Let $x\in\partial\Gamma$. Since $\partial \Gamma\setminus\{x\}$ is locally compact, we have the \emph{valency map} $\val\co\partial\Gamma\to\bN$ defined by
\[\val(x) := \#\text{ends}(\partial\Gamma\setminus\{x\}).\]
We define the set of \emph{local cut points} as
\[LC(\partial\Gamma):=\val^{-1}[2,\infty).\]

For two points $x,y\in LC(\partial\Gamma)$, we declare $x\sim y$  
if either  $x=y$ or 
\[
\val(x)=\val(y)=\#\pi_0(\partial\Gamma\setminus\{x,y\}).\]
In the latter case, we say $\{x,y\}$ is a \emph{cut pair};
since $\partial\Gamma$ has no global cut point we then have $x\ne y$.
It turns out that $\sim$ is an equivalence relation, the corresponding equivalence class of which will be denoted as $[x]$. 
So, it makes sense to define $\val[x]:=\val(x)$.

\begin{lem}\label{l:lc-gamma}
For each $x\in LC(\partial\Gamma)$, exactly one of the following holds.
\be[(i)]
\item  $\val(x)\in[3,\infty)$ and  $\#[x]=2$;
\item  $\val(x)=2$ and  $\#[x]=2$;
\item  $\val(x)=2$ and  $\#[x]=\infty$.
\ee
\end{lem}

\begin{rem}
In~\cite{Bowditch1998}, the alternative (i) is denoted as a \emph{$\approx$--pair}
and the symbol $\sim$ was reserved only for the other two alternatives.
\end{rem}

We now define
\[\Theta_1:=\{[x]\mid \#[x]=2\},\quad\Theta_2:=\{[x]\mid \#[x]=\infty\},\]
and set $T:=\Theta_1\cup \Theta_2$.
Let $A,B\sse\partial\Gamma$. We say $A$ \emph{separates} $B$, if $B$ intersects at least two distinct components of $\partial\Gamma\setminus A$.

Recall a Cantor set $\CC$ can be realized  as 
\[
\CC=S^1\setminus \coprod_{j\ge1} I_j\]
for some countable collection of open intervals $\{I_j\}$. 
This realization of $\CC$ is also called as a \emph{cyclically ordered Cantor set}.
Each pair of points $\partial I_j$ is called a \emph{jump} of the cyclically ordered Cantor set $\CC$.

\begin{prop}\label{p:sep}
Suppose $x\in LC(\partial\Gamma)$ satisfies $\val(x)=2$ and  $\#[x]=\infty$.
Then there exists a homeomorphism 
\[
h\co \overline{[x]}\to\CC\sse S^1\]
for some cyclically ordered Cantor set $\CC$
such that the following hold.
\be
\item
For some $y_1,y_2,\ldots\in LC(\partial\Gamma)$, we have that
each $h[y_i]$ is a jump and that
\[ \overline{[x]}\setminus[x] = \cup_{i\ge1} [y_i]
\sse\val^{-1}[3,\infty).\]
\item
For two $\theta,\xi\in T$, the class $\theta$ does not separate $\xi$.
\item
For four distinct points $x,y,z,w$ in $[x]$ whose images by $h$ appear cyclically in this order on $S^1$, we have that $\{x,z\}$ separate $\{y,w\}$ in $\partial\Gamma$.
\item
For $\theta,\xi,\eta\in T$, the class $\eta$ separate $\theta\cup\xi$ if and only if 
there exist two points $x,y\in\eta$ such that $\theta$ and $\xi$ are contained in distinct components of $\partial\Gamma\setminus\{x,y\}$.
\ee\end{prop}

We let $J(\overline{[x]})$ denote the set of all jumps in $\overline{[x]}$. 
Some of the jumps in $J(\overline{[x]})$ are missing from $[x]$ itself, while others are not.
Using the notations from the above proposition, we define 
\[
J_0(\overline{[x]}):=\{[y_i]\mid i\ge 1\}\]
as the set of ``missing jumps'' in $[x]$.

We say $F\sse T$ is \emph{inseparable} (or, \emph{null and full}) if there do not exist $\theta,\xi\in F$ and $\eta\in T$ such that $\eta$ separates $\theta\cup\xi$. A \emph{star} is a subset $F\sse T$ which is maximal among inseparable subsets.
We define $\Theta_3$ as the the set of stars with infinite cardinality.
For convention, we often identify a star with the union of the equivalence classes in it.
In particular,  each $\Theta_i$ is a collection of subsets of $\partial\Gamma$. 

Finally, we define $\Theta_1'$ as the collection of sets $\xi\sse\partial\Gamma$ satisfying both of the following properties:
\begin{itemize}
\item $\xi$ is a jump of some $\theta\in\Theta_2$ such that  $\xi\sse\theta$; in other words, $\xi\in J(\bar\theta)\setminus J_0(\bar\theta)$;
\item for some star $\eta\in\Theta_3$ containing $\theta$
and for the unique component 
 $U\in\pi_0(\partial\Gamma\setminus\xi)$ not intersecting $\theta$,
 we have that  $\eta \setminus\theta\sse U$.
\end{itemize}

For $A\sse\partial\Gamma$, we let $\Stab(A)$ denote the setwise stabilizer of $A$. 
Let
\[
\VV_1 = \{\Stab\theta\mid \theta\in \Theta_1\cup\Theta_1'\}.\]
Similarly, we let $\VV_i= \{\Stab\theta\mid \theta\in \Theta_i\}$ for $i=2,3$. 
We let $\VV=\VV_1\cup\VV_2\cup\VV_3$.
We note that $u\in\VV$ is elementary if and only if $u\in\VV_1$.

Recall an elementary splitting is uniquely determined by its collection of vertex stabilizer groups. The main result of~\cite{Bowditch1998} is the following.

\begin{thm}[Bowditch's canonical JSJ splitting]\label{thm:bowditch} 
Let $\Gamma$ be a one-ended nonelementary hyperbolic group that is not quasi-isometric to a cocompact Fuchsian group, and let $\VV$ be the collection of subgroups described above.
Then $\VV$ determines a maximal elementary splitting;
moreover, whenever $\Gamma$ splits over a two-ended group $H$, we can find some $v\in\VV_1\cup\VV_2$ such that $H\le v$.
\end{thm}

To describe this result more precisely, we define $\Sigma=(\VV,\EE)$ as a bipartite graph on $(\VV_1,\VV_2\cup \VV_3)$ such that the adjacency relation $\EE$ is defined as follows:
 $u\in\VV_1$ and $v\in\VV_2\cup\VV_3$ are adjacent if and only if $u\cap v$ is an infinite group (Proposition~\ref{prop:basic}).  We let $\Sigma$ be equipped with the natural conjugation action of $\Gamma$ defined as
\[
g.v:= gvg^{-1}.\]
Under this setting, Theorem~\ref{thm:bowditch} asserts that $\Sigma$ is a maximal elementary splitting. The following is now immediate.

\begin{cor} \label{cor:localcut}  If $\Gamma$ is hyperbolic and $\partial \Gamma$ contains a local cut point, then $\Gamma$ admits a splitting over a two-ended subgroup.\end{cor} 

See Haulmark \cite{Haulmarksplittings} for an extension of this to the relatively hyperbolic case.

\subsection{Algebraic and geometric features of the canonical splitting}
We have noted that every local cut point belongs to some cut pair. 
For a two-ended group $H\le \Gamma$, we define
\[
e(H):=\#\pi_0(\partial\Gamma\setminus\Lambda(H)).\]
For convention, we let $e(h):=e(\form{h})$ for a loxodromic $h$.
If $v\in\VV$, then we let $\delta(v)$ denote the set of incident edges on $v$
and $\deg(v):=\#\delta(v)$.
\begin{prop}[\cite{Bowditch1998}]
The following hold for a one-ended nonelementary hyperbolic group $\Gamma$ that is not quasi-isometric to a cocompact Fuchsian group
and for its canonical JSJ splitting $\Sigma=(\VV=\VV_1\cup\VV_2\cup\VV_3,\EE)$.
\be
\item
Every local cut point $x\in\partial\Gamma$ is contained in the limit set of some $u\in\VV_1\cup\VV_2$;
furthermore, if $\#[x]=2$, then $u$ can be chosen to be 2-ended.
\item 
If $u\in\VV_1$, then we have $e(u)=\val(u)=\deg(u)\in[2,\infty)$;
furthermore, if $\deg(u)=2$, then some vertex of $\VV_3$ is adjacent to $u$.
\item
If $u\in\VV_2\cup\VV_3$, then $\deg(u)=\infty$.
\item 
Every loxodromic element $\gamma\in\Gamma$ with $e(\gamma)>1$ 
belongs to some $u\in \VV_1\cup \VV_2$ such that $\val(u)=e(\gamma)$.
\ee
\end{prop}

\begin{defn} \label{def:maxhang} Each group $u\in\VV_2$ is called a \emph{maximal hanging Fuchsian (MHF)} group. \end{defn} 
For an MHF group $u\in\VV_2$, one can find a discrete representation with finite kernel
\[\rho\co u\to\Isom(\mathbb{H}^2)\]
and an equivariant cyclic-order-preserving homemorphism \[h\co \Lambda(u)\to \Lambda(\rho(u))\sse \partial \mathbb{H}^2.\]
In particular, $\rho(u)$ is naturally identified with the orbifold fundamental group of a (not necessarily oriented) compact hyperbolic two--orbifold $S$.
We have a finite disjoint union 
\[
\partial S = \coprod_{i=1}^m  \partial_i S,\]
where each component $\partial_i S$ is either a circle $S^1=\bR/\bZ$ or a compact interval with mirrors $\bR/(\bZ_2\ast\bZ_2)$.
Let us denote by $P_i\le u$ the preimage of $[\partial_i S]\in\pi_1^{\mathrm{orb}}(S)=\rho(u)$, with an arbitrary choice of the base point. 
Then there exists a natural one-to-one correspondence between the incident edges on $u$
and the $u$--conjugates of all $P_i$, which maps an edge $e$ to its stabilizer group $\Gamma(e)$.
In particular, $\deg(u)=\infty$.

We call the collection of groups $\{P_i^g\mid g\in u\text{ and }i=1,2,\ldots,m\}$ as the \emph{peripheral structure} of the MHF group $u$. Each element of $P_i^g$ is called \emph{peripheral}.
 This algebraic feature of $u\in\VV_2$, along with quasi-convexity and fullness, actually characterizes the maximal hanging Fuchsian groups.
We note that if the limit sets of two distinct MHF groups $u,v\in\VV_2$ are not disjoint, then they share a common missing jump in $J_0(\Lambda u)\cap J_0(\Lambda v)$. In particular, the midpoint of the interval $[u,v]$ in $\Sigma$ has degree at least three.

\bd \label{d:rigidrel} 
Let $\Sigma=(\VV,\EE)$ be an elementary splitting of $\Gamma$. We say a vertex $v$ of  $\Sigma$ is \emph{rigid (rel incident edges)} if there does not exist a graph of groups decomposition $\GG$ of 
 $\Gamma(v)$ in such a way that each group in the set \[\{ \Gamma(e)\mid e\in \delta(v)\}\] is $\Gamma(v)$--conjugate into a  vertex stabilizer group in $\Sigma'$
and that each edge group of $\GG$ is two-ended.
\ed
\begin{rem}We note that $\Gamma(v)$ is not required to be one-ended.\end{rem}

Each vertex $u\in \VV_3$ is a quasi-convex, full, nonelementary, non-MHF group such that it is rigid rel  incident edges in the canonical splitting of $\Gamma$. This algebraically characterizes the groups in $\VV_3$. One can also see that a rigid vertex group is rigid relative to the collection of incident two-ended edge groups, in the sense of Definition~\ref{d:rigid}.

For a two-ended subgroup $H\le\Gamma$,
the value $e(H)$ coincides with the maximum number of ends of the group pair $(\Gamma,H')$ where $H'$ ranges over finite index subgroups of $H$. In this context, one can realize $\VV_1$ as the collection of maximal elementary groups containing loxodromics $\gamma$ such that the number of ends of a group pair $(\Gamma,\form{\gamma})$ is larger than 1
and such that $\gamma$ is not contained in a MHF group as a non-peripheral element.

\section{Concrete examples of the Planarity conjectures} 
In this section, we illustrate the validity of planarity conjectures for doubles and limit groups. These results follow from (among other places) Haissinsky's work \cite{Haissinsky2015IM} on planar boundaries.  The most relevant result to our discussion here is 

\begin{thm} \label{thm:Hnosier} \cite[Corollary 1.14]{Haissinsky2015IM} Every hyperbolic group $G$ with a one-dimensional planar boundary and no elements of order two is virtually Kleinian if and only if every carpet group is virtually Kleinian. In particular, if $G$ has no carpet subgroup, then G is virtually Kleinian.\end{thm} 

A {\it carpet subgroup} is a quasi-convex subgroup $H$ of a hyperbolic group such that $\partial H $ is a Sierpinski carpet. 
A {\it Sierpinski carpet} is a planar 1-dimensional  Peano continuum without local or global cut points.  The Sierpinski carpet occurs as the boundary of a hyperbolic group that is the fundamental group of a hyperbolic 3-manifold with totally geodesic boundary.  If a hyperbolic group $H$ has boundary a Sierpinski carpet, then it does not split over a finite or a cyclic group by our previous discussion.   A hyperbolic group with no elements of order two admits a finite hierarchy over cyclic and finite groups \cite{loudertouikan}.  

Hyperbolic doubles of free groups and hyperbolic limit groups are examples of hyperbolic groups which do not contain a carpet subgroup.  They admit a hierarchy over finite and cyclic groups which terminates in free groups. This hierarchy is particularly simple for doubles. The study of relatively hyperbolic boundaries of free groups, where some words in the free group are parabolic, was first initiated by Otal, although he did not use this language. 
 
\subsection{Otal's results}\label{sec:Otal} 
The main result of Otal \cite[Theorem 1]{Otal1992} concerns collections of conjugacy classes in a free group $F$. 
Let $P = \pbr{ \gamma_1,...,\gamma_n }$ be a multiword (i.e. a collection of conjugacy classes of root-free words) such that no nontrivial powers of two such words are conjugate to each other. 
Note that this implies the collection of infinite cyclic subgroups generated by $\gamma_1,...\gamma_n$ form an almost malnormal collection,
as defined in Section~\ref{sec:prelim}.
We then let $\mathcal{P}$ be the collection of conjugates of words in $P$.  Then the pair $(F, \mathcal{P})$ is a relatively hyperbolic group and the Bowditch boundary $\partial_B(F, \mathcal{P})$ is a quotient of the Cantor set $\partial F$  obtained by identifying the endpoints of conjugates of words in $P$.
Otal denoted this Bowditch boundary $\partial_B(F,\mathcal{P})$ as $K_P$, which Cashen also called the \emph{decomposition space} of $(F, \mathcal{P})$~\cite{Cashen2016AGT}.

We will also need the concept of a {\it relative splitting of a relatively hyperbolic group pair $(G, \mathcal{P})$} .  This is a splitting of $G$ as a free product with amalgamation  $A *_C B$ or HNN extension $A *_C$ such that every $P \in \mathcal{P}$ is in $A$ or $B$. 
Putting Otal's results in the language of relatively hyperbolic groups we have:

\begin{thm}[\cite{Otal1992,Cashen2016AGT}]
Let $F$ be free of rank at least 2 and $(F, \mathcal{P})$ a relatively hyperbolic group pair where each $P \in \mathcal{P}$ is infinite cyclic. Suppose that $(F, \mathcal{P})$ does not admit a relative splitting over a virtually cyclic group.  Then if $\partial_B(F, \mathcal{P})$ is planar, $F$  is the fundamental group of a handlebody $H$ where the conjugacy classes of $\mathcal{P}$ correspond to a multicurve which is isotopic into $\partial H$. 
 \end{thm} 

In other words, the collection of conjugacy classes $\mathcal{P}$ is geometric; see Section~\ref{ex:double}.


\subsection{Doubles and limit groups}
Let us now consider an indecomposible, acylindrical, root-free word $w$ in $F$, such that the (hyperbolic) group boundary $\partial D(F,w)$ is planar. We will show that planarity of the boundary, along with Bowditch's and Otal's results,  implies that the double $D(F,w)$ is convex cocompact Kleinian.

We let $\PP$ denote the set of conjugates of  $\form{w}$. We have seen by work of Tran that the Bowditch boundary $\partial_B(F, \PP)$ is obtained from $\partial F$ by identifying the pair $gw^\infty g^{-1}$ with $gw^{-\infty} g^{-1}$ for each $g\in F$. Fix an embedding $\partial D(F,w)\hookrightarrow S^2$.
Then there exists a $D(F,w)$--equivariant collection of arcs $\gamma_g$ joining the pair 
\[
Z_g:=\left(gw^\infty g^{-1},gw^{-\infty} g^{-1}\right)\]
for each $g\in F$ in such a way that $\gamma_g$ belongs to the component of $\partial D(F,w)\setminus Z_g$ not containing $\partial F$. Indeed, $F$ is the stabilizer of a vertex in the maximal splitting defined by Bowdtich described in Section 4, and we denote this stabilizer by $F_v$. The pairs $Z_g$ are the edge groups associated to edges emanating from $v$. Each of these pairs is an equivalence class of local cut points and each pair separates the boundary into two (path connected) components, with the limit set of $F_v$ contained in one component.  Each of these cut pairs does not separate any of the other cut pairs by work of Bowditch, see (1) of Proposition~\ref{p:sep}.  Therefore we can connect each pair $Z_g$ in the path component not containing $F_v$ by a path which we call $\gamma_g$.  The other $Z_g$ pairs do not meet this path component, so the collection of arcs is disjoint.  
Contracting $\gamma_g$ to a point for each $g$, one sees from a classical result of Moore that $\partial_B(F, \PP)$ is planar as realized in Lemma~\ref{lem:Tranlemma}. 
See~\cite{Otal1992} and Cashen ~\cite{Cashen2016AGT} for more details on similar arguments. 

Once we see that $\partial_B(F, \PP)$ is planar, we deduce from Otal's result that $w$ is geometric. We have seen in Section~\ref{ex:double} that the geometricity implies that $D(F,w)$ is actually a convex cocompact Kleinian group.

Hyperbolic doubles of free groups are special cases of  limit groups.  
Recall a finitely generated group $L$ is called a \emph{limit group} (or, a \emph{fully residually free group}) if for each finite subset $A\sse L$ there exists a homomorphism
\[
\phi_A\co L\to F\]
 to a fixed nonabelian free group $F$ such that the restriction of $\phi_A$ to $A$ is injective.

A torsion--free finitely generated group $G$ is said to admit a \emph{cyclic hierarchy (of level at most $d$ over free groups)}  if one of the following conditions are satisfied. \be \item $d=0$ and $G$ is free; \item $d>0$ and $G$ splits as a finite graph of vertex groups $\{G_i\}$ with  cyclic (possibly trivial) edge groups, so that each $G_i$ admits a cyclic hierarchy over free groups of level at most $d-1$. \ee In particular, the double of a free group admits a cyclic hierarchy of level $1$. 
More generally, all hyperbolic limit groups admit cyclic hierarchies over free groups;  in fact, they are virtually free-by-cyclic~\cite{HW2010GD}.

Let $L$ be a hyperbolic limit group. 
Then every nontrivial quasi-convex subgroup admits a cyclic hierarchy as well.
On the other hand, Sierpinski carpet groups do not split over cyclic groups (including trivial) since their boundaries are connected and do not have local cut points. 
In particular, $L$ does not contain a carpet group.
Haissinsky's work (Theorem  \ref{thm:Hnosier} here) implies that hyperbolic limit groups with planar boundary are virtually Kleinian. 

Even when a hyperbolic group is torsion free, and hence acts effectively on its boundary, it may be virtually Kleinain without being Kleinian.  This can happen if the action does not extend to the whole of $S^2$.  Indeed, Kapovich and Kleiner gave an example of such phenomenon in \cite[Section 8]{KK2000}. See \cite{HSTAJM} for more examples.  All of these examples have boundaries which split over a two-ended group, and we remark that this condition is necessary. A virtually Kleinian group is quasi-isometric to a Kleinian group.

\begin{prop} \label{prop:virtual} Let $G$ be a torsion-free hyperbolic group that is quasi-isometric to a Kleinian group.  Then if $G$ does not split over a finite or 2-ended group, $G$ can be realized as a Kleinian group. \end{prop} 

Suppose that $G$ is a torsion-free hyperbolic group that is quasi-isometric to a Kleinian group.
Note that $G$ acts effectively on its boundary since it is torsion free. 
Since any Kleinian group that is hyperbolic can be realized as a convex-cocompact Kleinian group, its boundary is planar.  If $\partial G \simeq S^1$, then $G$ is virtually Fuchsian, by work of Tukia \cite{Tukia1988}, Gabai \cite{Gabaiconvergence} and Casson--Jungreis \cite{CassonJungreis}. 
In this case $G$ can be realized as a Fuchsian subgroup of $\Isom(\mathbb{H}^2) \leq \Isom(\mathbb{H}^3)$. 

If $\partial G$ is not $S^1$, and $G$ does not split over a 2-ended group, then $\partial G$ does not have any local cut points by Bowditch's work (Corollary \ref{cor:localcut}).
We now have that $\partial G$ is a planar Peano continuum without local or global cut points. 
Let us regard $\partial G$ as a subspace of $S^2$. 
Assume first that $\dim \partial G = 2$. Then since $\partial G$ is a subset of a 2-dimensional manifold, $\partial G$ contains an open 2--disk; see \cite[Corollary 1, page 46]{HWdim}.
Since $G$ acts on $\partial G$ with dense orbits and by homeomorphisms, this implies that $\partial G$ is open in $S^2$. Since $\partial G$ is compact, we see that $\partial G = S^2$. 
We have assumed that  $G$ is torsion-free and quasi-isometric to a Kleinian group.  It is a result of Cannon and Cooper, \cite{CannonCooper}, using work of Sullivan \cite{Sullivan}, that $G$ acts geometrically on $\mathbb{H}^3$. 

Finally, suppose that $\partial G$ is 1-dimensional (see \cite[Theorem 4]{KK2000} for a classification in the 1-dimensional case). As we have noted, the boundary of $G$ is a Sierpinski carpet in this case.
The double of $G$ along the subgroups that stabilize the non-separating circles is a hyperbolic group $\widehat G$ with boundary $S^2$.  See \cite[Theorem 5 and Section 5]{KK2000}.
Since $G$ is virtually Kleinian, so is $\widehat G$. 
Therefore, by the above argument, $\widehat G$ can be realized as a Kleinian group.  It follows that $G$ is also Kleinian. 


 

\section{Acknowledgements} 
We thank the Centre International de Rencontres Math\'ematiques (CIRM), where these discussions were started when both authors were visiting the institute. 
The first named author was supported by the National Research Foundation funded by the government of Korea (2018R1A2B6004003).
The second  named author was partially supported under NSF DMS- 1709964.







\begin{thebibliography}{10}

\bibitem{BF1992}
M.~Bestvina and M.~Feighn, \emph{A combination theorem for negatively curved
  groups}, J. Differential Geom. \textbf{35} (1992), no.~1, 85--101.

\bibitem{BM1991}
M.~Bestvina and G.~Mess, \emph{The boundary of negatively curved groups}, J.
  Amer. Math. Soc. \textbf{4} (1991), no.~3, 469--481. \MR{1096169}

\bibitem{BowditchGF}
B.~H. Bowditch, \emph{Geometrical finiteness for hyperbolic groups}, J. Funct.
  Anal. \textbf{113} (1993), no.~2, 245--317. \MR{1218098 (94e:57016)}

\bibitem{Bowditch1998}
\bysame, \emph{Cut points and canonical splittings of hyperbolic groups}, Acta
  Math. \textbf{180} (1998), no.~2, 145--186.

\bibitem{Bowconnected}
\bysame, \emph{Connectedness properties of limit sets}, Trans. Amer. Math. Soc.
  \textbf{351} (1999), no.~9, 3673--3686. \MR{1624089 (2000d:20056)}

\bibitem{Bowditchrelhyp}
\bysame, \emph{Relatively hyperbolic groups}, Internat. J. Algebra Comput.
  \textbf{22} (2012), no.~3, 1250016, 66. \MR{2922380}

\bibitem{BH1999}
M.~R. Bridson and A.~Haefliger, \emph{Metric spaces of non-positive curvature},
  Grundlehren der Mathematischen Wissenschaften [Fundamental Principles of
  Mathematical Sciences], vol. 319, Springer-Verlag, Berlin, 1999.

\bibitem{Canary1993}
R.~D. Canary, \emph{Ends of hyperbolic {$3$}-manifolds}, J. Amer. Math. Soc.
  \textbf{6} (1993), no.~1, 1--35.
  
  \bibitem{CannonCooper}
  J. ~W. Cannon and D. Cooper, \emph{A characterization of cocompact hyperbolic and finite-volume hyperbolic groups in dimension three}, 
Trans. Amer. Math. Soc.\textbf{330}, no. ~1, 419--431.

\bibitem{Cashen2016AGT}
C.~r~H. Cashen, \emph{Splitting line patterns in free groups}, Algebr. Geom.
  Topol. \textbf{16} (2016), no.~2, 621--673. \MR{3493403}

\bibitem{CassonJungreis}
A.~J. Casson and D.~Jungreis, \emph{Convergence groups and Seifert fibered
  3-manifolds}, Invent. Math \textbf{118} (1994), 441--456.

\bibitem{Gabaiconvergence}
D.~Gabai, \emph{Convergence groups are fuchsian groups}, Ann. Math.
  \textbf{136} (1992), 447--510.

\bibitem{Floydsquash}
V.~Gerasimov and L.~Potyagailo, \emph{Quasi-isometric maps and floyd boundaries
  of relatively hyperbolic groups}, Journal of the European Mathematical
  Society \textbf{15} (2009), no.~6.

\bibitem{GW2010}
C.~McA. Gordon and H.~Wilton, \emph{On surface subgroups of doubles of free
  groups}, J. Lond. Math. Soc. (2) \textbf{82} (2010), no.~1, 17--31.
  \MR{2669638 (2011k:20085)}

\bibitem{Gromovhyp}
M.~Gromov, \emph{Hyperbolic groups}, Essays in group theory, Math. Sci. Res.
  Inst. Publ., vol.~8, Springer, New York, 1987, pp.~75--263.

\bibitem{GMS}
D.~ Groves, J.~F~ Manning, and A. Sisto, \emph{Boundaries of
  {D}ehn fillings}, Geometry and Topology (to appear).

\bibitem{HW2010GD}
M. Hagen and D.l~T. Wise, \emph{Special groups with an elementary
  hierarchy are virtually free-by-{$\Bbb Z$}}, Groups Geom. Dyn. \textbf{4}
  (2010), no.~3, 597--603. \MR{2653976}

\bibitem{Haissinsky2015IM}
P. Ha\"{i}ssinsky, \emph{Hyperbolic groups with planar boundaries}, Invent.
  Math. \textbf{201} (2015), no.~1, 239--307. \MR{3359053}

\bibitem{Haulmarksplittings}
M.~ {Haulmark}, \emph{{Local cut points and splittings of relatively
  hyperbolic groups}}, arXiv e-prints (2017), arXiv:1708.02855.

\bibitem{Burns}
B.~B.~ Healy, \emph{Rigidity properties for hyperbolic generalizations},
   (2018), Preprint, arXiv:1803.10153.

\bibitem{HruskaRH}
G.~C. Hruska, \emph{Relative hyperbolicity and relative quasiconvexity for
  countable groups}, Algebr. Geom. Topol. \textbf{10} (2010), no.~3,
  1807--1856. \MR{2684983 (2011k:20086)}

\bibitem{HWpreprint}
G.~C. Hruska and G.~Walsh, \emph{Relatively hyperbolic groups with planar
  boundary},  (2019), Preprint.

\bibitem{HSTAJM}
G.~C. Hruska, E. Stark, and H.~C. Tran, \emph{Surface group amalgams that
  (don't) act on 3-manifolds}, American Journal of Mathematics, to appear,
  arXiv:1705.01361.

\bibitem{HWdim}
W. Hurewicz and H. Wallman, \emph{Dimension theory}, Princeton
  University Press, Princeton, N. J., 1941.

\bibitem{KapBenakli}
I. Kapovich and N. Benakli, \emph{Boundaries of hyperbolic groups},
  Combinatorial and geometric group theory ({N}ew {Y}ork, 2000/{H}oboken, {NJ},
  2001), Contemp. Math., vol. 296, Amer. Math. Soc., Providence, RI, 2002,
  pp.~39--93. \MR{1921706}

\bibitem{KK2000}
M.~Kapovich and B.~Kleiner, \emph{Hyperbolic groups with low-dimensional
  boundary}, Ann. Sci. \'{E}cole Norm. Sup. (4) \textbf{33} (2000), no.~5,
  647--669. \MR{1834498}

\bibitem{Kapovicharticle}
M.~Kapovich \emph{Lectures on quasi-isometric rigidity},  PCMI Lecture Note Series, (2012).

\bibitem{KapKleinlow}
M.~Kapovich and B.~Kleiner, \emph{Hyperbolic groups with
  low-dimensional boundary}, Ann. Sci. \'{E}cole Norm. Sup. (4) \textbf{33}
  (2000), no.~5, 647--669. \MR{1834498}

\bibitem{Kapprob}
M.~ Kapovich, \emph{Problems on boundaries of groups and kleinian groups},
  2007.

\bibitem{KW2012}
{S.-h.} Kim and H.~Wilton, \emph{Polygonal words in free groups}, Q. J.
  Math. \textbf{63} (2012), no.~2, 399--421. \MR{2925298}

\bibitem{loudertouikan}
L. Louder and N.~ Touikan, \emph{Strong accessibility for finitely
  presented groups}, Geometry and Topology \textbf{21}, 1805?1835.

\bibitem{Manningsquash}
J.~F. Manning, \emph{The bowditch boundary of $({G},\mathcal{H})$ when ${G}$ is
  hyperbolic}, arXiv:1504.03630.

\bibitem{Morgan}
J.~W. Morgan, \emph{On {T}hurston's uniformization theorem for
  three-dimensional manifolds}, The {S}mith conjecture ({N}ew {Y}ork, 1979),
  Pure Appl. Math., vol. 112, Academic Press, Orlando, FL, 1984, pp.~37--125.
  \MR{758464}

\bibitem{Otal1992}
J.-P. Otal, \emph{Certaines relations d'\'{e}quivalence sur l'ensemble
  des bouts d'un groupe libre}, J. London Math. Soc. (2) \textbf{46} (1992),
  no.~1, 123--139. \MR{1180888}

\bibitem{Scott1979}
P. Scott and T. Wall, \emph{Topological methods in group theory},
  Homological group theory ({P}roc. {S}ympos., {D}urham, 1977), London Math.
  Soc. Lecture Note Ser., vol.~36, Cambridge Univ. Press, Cambridge, 1979,
  pp.~137--203. \MR{564422 (81m:57002)}

\bibitem{seriescrash}
C. Series, \emph{A crash course on {K}leinian groups}, Rend. Istit. Mat.
  Univ. Trieste \textbf{37} (2005), no.~1-2, 1--38 (2006). \MR{2227047}

\bibitem{Serre-trees}
J.-P Serre, \emph{Trees}, Springer Monographs in Mathematics,
  Springer-Verlag, Berlin, 2003, Translated from the French original by John
  Stillwell, Corrected 2nd printing of the 1980 English translation.
  \MR{1954121}

\bibitem{Short1990}
H. Short, \emph{Groups and combings},  (1990).


\bibitem{Sullivan}
D. Sullivan, \emph{On the ergodic theory at infinity of an arbitrary discrete group of hyperbolic
motions}, Proc. Stony Brook Conf., Ann. of Math. Studies No. 97, Princeton Univ. Press,
Princeton, N.J., 1981, pp. 465-496.


\bibitem{Transquash}
H.~C. Tran, \emph{Relations between various boundaries of relatively
  hyperbolic groups}, Internat. J. Algebra Comput. \textbf{23} (2013), no.~7,
  1551--1572. \MR{3143594}

\bibitem{Tukia1988}
P.~Tukia, \emph{Homeomorphic conjugates of fuchsian groups}, J. Reine Angew.
  Math. \textbf{391} (1988), 1--54.


\end{thebibliography}

\providecommand{\bysame}{\leavevmode\hbox to3em{\hrulefill}\thinspace}
\providecommand{\MR}{\relax\ifhmode\unskip\space\fi MR }
\providecommand{\MRhref}[2]{%
  \href{http://www.ams.org/mathscinet-getitem?mr=#1}{#2}
}
\providecommand{\href}[2]{#2}

\end{document}